\documentclass[twoside,11pt]{article}
\usepackage{graphicx,amsmath,amsthm,amssymb,latexsym,amsfonts,color}
\usepackage{subfigure}
\usepackage[bookmarksnumbered, colorlinks, plainpages]{hyperref}
\usepackage{lscape}
\newtheorem{thm}{\bf Theorem}
\newtheorem{lem}{\bf Lemma}

\def\h{\hspace{-0.2cm}}
\def\dis{\displaystyle}

\newcommand{\bx}{{\bf x}}

\begin{document}
\title{A  new preconditioner for elliptic PDE-constrained optimization problems}

\author{{ Hamid Mirchi and Davod Khojasteh Salkuyeh}\\[2mm]
\textit{{\small Faculty of Mathematical Sciences, University of Guilan, Rasht, Iran}} \\
\textit{{\small E-mails: hamidmirchi@gmail.com and khojasteh@guilan.ac.ir}\textit{}}}
\date{}
\maketitle
\noindent{\bf Abstract.}  We propose a  preconditioner to accelerate the convergence of the GMRES iterative method  for solving the system of linear equations obtained from  discretize-then-optimize approach applied to optimal control problems constrained by a partial differential equation.
Eigenvalue distribution of the preconditioned matrix as well as its eigenvectors are discussed. Numerical results of the proposed preconditioner are compared with several existing preconditioners to show its efficiency.\\[-3mm]

\pagestyle{myheadings}\markboth{H. Mirchi, D.K. Salkuyeh}{A new preconditioner for elliptic PDE-constrained optimization problems}

\thispagestyle{empty}

\section{Introduction} \label{SEC1}

Consider the following linear elliptic distributed  optimal control problem:
\begin{eqnarray}
\min_{u,f}~~~ &\h\h&\frac{1}{2} \|u-u_*\|_{L^2(\Omega)}^2+\beta\|f\|_{L^2(\Omega)}^2,\label{Eq001}\\
\text{s.t.}~~~~&\h\h& -\Delta u=f\quad \text{in}~~ \Omega,\label{Eq002}\\
&\h\h&u=g \quad \text{on}~~\partial\Omega_1\quad \text{and}~~~ \frac{\partial u}{\partial n}=g\quad \text{on}~~ \partial\Omega_2,\label{Eq0003}
\end{eqnarray}
where $\partial\Omega$ is the boundary of  domain $\Omega$ in $\Bbb{R}^2$ or $\Bbb{R}^3$, $\partial \Omega_1  \cup  \partial\Omega_2=\partial\Omega$ and $\partial \Omega_1  \cap \partial \Omega_2=\phi$. Such problems, introduced by Lions in \cite{Lions1968}, consists in a cost functional  \eqref{Eq001} to be minimized subject to a partial differential problem \eqref{Eq002}-\eqref{Eq0003} in the domain $\Omega$.
Here, the  ``desired state'' $u_*$ is a known function, and we need to find $u$ which satisfies the PDE problem that it is close to $u_*$ as possible in the $L_2$-norm.  Furthermore, $\beta>0$ is a regularization parameter.

There are two options for solving the problem: discretize-then-optimize and optimize-then-discretize. Using the discretize-then-optimize approach, Rees et al. in \cite{Rees2010},  transformed the problem into a linear system of saddle-point form. They applied the Galerkin finite-element method to the weak formulation of Eqs. \eqref{Eq001}-\eqref{Eq0003} and obtained the finite-dimensional optimization problem
\begin{equation}\label{Eq004}
\left\{
\begin{array}{l}
\dis\min_{{\bf u},{\bf f}}~~~ \frac{1}{2}{\bf u}^TM{\bf u}-{\bf u}^T{\bf b}+\|{\bf u}_*\|_2^2+\beta {\bf f}^TM{\bf f}, \\[2mm]
\text{s.t.}~~~~ K{\bf u}=M{\bf f}+{\bf d},
\end{array}
\right.
\end{equation}
where $K,M\in\Bbb{R}^{m\times m}$ are the stiffnes (the discrete Laplacian) and mass matrices,  ${\bf d}\in\Bbb{R}^{m}$ contains the terms coming from the boundary values of the discrete solution, and ${\bf b}\in\Bbb{R}^{m}$ is the Galerkin projection of the discrete state $u_*$. Applying the Lagrange multiplier technique to the minimization problem \eqref{Eq004}, and selecting the $({\bf f},{\bf u},\lambda)$ ordering of the unknowns, results in the following linear system of equations
\begin{equation}\label{MainSys}
\mathcal{A}{\bf x} =
\begin{pmatrix}
2\beta M & 0 & -M \\
0        & M & K^T \\
-M       & K &  0
\end{pmatrix}
\begin{pmatrix}
{\bf f} \\
{\bf u}\\
\lambda
\end{pmatrix}
=
\begin{pmatrix}
0 \\
{\bf b}\\
{\bf d}
\end{pmatrix}={\bf g},
\end{equation}
where $\lambda\in\Bbb{R}^m$ is a vector of Lagrange multipliers.

Both of matrices $M$ and $K$ are sparse and hence the matrix $\mathcal{A}$ is sparse. Further, the matrix $M$ is symmetric positive definite. In general the matrix $\mathcal{A}$ is symmetric and indefinite. Therefore,  the MINRES iterative method presented by Paige and Saunders in \cite{MINRES}   can  be  used for solving the system  \eqref{MainSys}. Several preconditioners have been presented for solving the linear system \eqref{MainSys} in the literature.

In \cite{Rees2010}, Rees et al. applied the MINRES  algorithm  in conjunction with the block diagonal preconditioner
\[
 \mathcal{P}_{D}= \begin{pmatrix}
2\beta M & 0 & 0 \\
0       & M & 0 \\
0       & 0 &  KM^{-1}K^T
\end{pmatrix}.
\]
The projected preconditioned  conjugate gradient (CG) method combined with the constraint preconditioner
\[
\mathcal{P}_{C}= \begin{pmatrix}
0 & 0 & -M \\
0 & 2\beta K^TM^{-1}KM & K^T \\
-M & K &  0
\end{pmatrix},
\]
was presented by Gould et al. in \cite{PPCG}.  Rees and Stoll in \cite{Rees-NLWA} proposed the block-triangular (BT) preconditioner
\[
\mathcal{P}_{BT}= \begin{pmatrix}
	2\beta M & 0 & 0 \\
	0 & M & 0 \\
	-M & K &  KM^{-1}K^T
\end{pmatrix}.
\]
Bai in \cite{Bai-BCD}, presentated the GMRES methods coupled with the  block-counter-diagonal (BCD)  preconditioner
\begin{equation*}
\mathcal{P}_{BCD}= \begin{pmatrix}
	0 & 0 & -M \\
	0 & M & 0 \\
	-M & 0 &  0
\end{pmatrix},\\
\end{equation*}
and block-counter-tridiagonal (BCT) preconditioner
\[
\mathcal{P}_{BCT}= \begin{pmatrix}
	0 & 0 & -M \\
	0 & M & K^T \\
	-M & K &  0
\end{pmatrix}.\\
\]
 Zhang and Zheng in \cite{Zhang-JCM} introduced the block-symmetric (BS)  preconditioner
\[
\mathcal{P}_{BS} = \begin{pmatrix}
	2\beta M & 0 & -M \\
	0 & M & 0 \\
	-M & 0 &  0
\end{pmatrix},\\
\]
and the block-lower-triangular preconditioner
\[
\mathcal{P}_{BLT} = \begin{pmatrix}
	2\beta M & 0 & 0 \\
	0 & M & 0 \\
	-M & K &  -\frac{1}{2\beta}M
\end{pmatrix}.
\]
Recently, Ke and Ma in \cite{Ke-CAMWA}, proposed the following four preconditioners
\begin{eqnarray*}
\mathcal{P}_{1}&\h=\h& \begin{pmatrix}
2\beta M & 0 & -M \\
0 & 0    & K^T \\
-M & K &  0
\end{pmatrix},\quad \mathcal{P}_{2}=\begin{pmatrix}
2\beta M & 0 & -M \\
0 & M & K^T \\
0 & K &  0
\end{pmatrix},\\
\mathcal{P}_{3}&\h=\h& \begin{pmatrix}
2\beta M & 0 & -M \\
0 & M    & 0 \\
-M & K &  0
\end{pmatrix},\quad \mathcal{P}_{4}= \begin{pmatrix}
2\beta M & 0 & -M \\
0 & M & K^T \\
-M & 0 &  0
\end{pmatrix}.
\end{eqnarray*}
The preconditioners $\mathcal{P}_{3}$ and $\mathcal{P}_{4}$ are suitable for regularization parameters being very small ($\beta< 10^{-6}$), and the preconditioners $\mathcal{P}_{1}$ and $\mathcal{P}_{2}$ for $\beta\geqslant 10^{-6}$.

In the implementation of a preconditioner $\mathcal{M}$ in a Krylov subspace method like GMRES for the system \eqref{MainSys}, we should compute vectors of the form $z=\mathcal{M}^{-1}r$,  where  $r\in\mathbb{R}^{n}$. This  can be done by solving the system $\mathcal{M}z=r$.  If $\mathcal{M}=\mathcal{P}_D$, then for solving the system $\mathcal{P}_Dz=r$ we should solve three sub-systems with the coefficient matrix $M$.  In the case that the matrix $M$ is symmetric positive definite, these systems can be solved exactly using the Cholesky factorization of  $M$ or inexactly using the CG method or its preconditioned version (PCG). Similar strategy can be implemented for other preconditioners. For all the reviewed  preconditioners the coefficient matrices of the sub-systems are $M$, $K$ or $K^T$. Since, in our numerical results the matrices $M$ and $K$ are symmetric positive definite, the same strategy as the preconditioner $\mathcal{P}_D$ can be used for solving these systems.

In this paper, we present a new  preconditioner for solving  the system \eqref{MainSys} and investigate the eigenvalue distribution of the preconditioned system.
Numerical comparison  with the recently presented preconditioners are presented.

This paper is organized  as follows. In Section \ref{Sec2} we  present the new preconditioner and  investigate its  properties.  Section \ref{Sec3} is devoted to some numerical results. Concluding remarks are  presented in  Section \ref{Sec4}.

\section{The new preconditioner} \label{Sec2}

Now, we propose the preconditioner
\begin{equation}\label{NonsinP}
\mathcal{P}=\begin{pmatrix}
0 &  K & 0 \\
0      &  M & K^T \\
-M      & K & 0
\end{pmatrix},
\end{equation}
for the system \eqref{MainSys}. In fact we want to solve the preconditioned system $\mathcal{P}^{-1}\mathcal{A}{\bf x}=\mathcal{P}^{-1}{\bf g}$ by a Krylov subspace method like GMRES.  In general, favorable convergence rates of Krylov subspace methods like GMRES are often obtained by a clustering of most of the eigenvalues of preconditioned matrices around $(1, 0)$ and away from zero \cite{BenziSurvey}. In the sequel we investigate the nonsingularity of the matrix $\mathcal{P}$ and the spectrum of the preconditioned matrix $\mathcal{P}^{-1}{\mathcal{A}}$.

\begin{lem}
If the matrices $M$ and  $K$  are nonsingular, then the matrix  $\mathcal{P}$ defined in \eqref{NonsinP} is nonsingular.
\begin{proof}
	Let ${\bf w}=(x;y;z)\in \Bbb{R}^{3m}$ and $\mathcal{P}{\bf w}=0$. In this case, we have
	\begin{eqnarray}
	Ky&\h=\h&0, \label{Nonsin1}\\
	My+K^Tz&\h=\h&0,\label{Nonsin2}\\
	-Mx+Ky&\h=\h&0.\label{Nonsin3}
	\end{eqnarray}
Since the matrix $K$ is nonsingular, from Eq.  \eqref{Nonsin1} we get $y=0$. Substituting $y=0$ into Eqs. \eqref{Nonsin2} and
\eqref{Nonsin3}, from the nonsingularity of the matrices  $K$ and $M$, we see that $x=z=0$.  Therefore,  ${\bf w}=0$ which completes the proof.
\end{proof}
\end{lem}

\begin{thm}\label{Thm1}
Assume that the matrices $M$ and $K$ are nonsingular.  Then, the matrix $\mathcal{P}^{-1}{\mathcal{A}}$ has an eigenvalues 1 with algebraic multiplicity $2m$ and the corresponding eigenvectors are of the form
\begin{equation}\label{Eigw1}
{\bf w}=\begin{pmatrix}
\frac{1}{2\beta}(z+M^{-1}Ky) \\
	y  \\
	z
\end{pmatrix},
\end{equation}
for any $y,z\in\Bbb{C}^{m}$ with at least one of them being nonzero. The remaining eigenvalues are of the form
\begin{equation}\label{Eignonunit}
\lambda=2\beta+\frac{z^*Mz}{z^*KM^{-1}K^Tz},
\end{equation}
and the corresponding eigenvectors are of the form
\[
v=\begin{pmatrix}
-M^{-1}KM^{-1}K^Tz\\
-M^{-1}K^Tz  \\
z
\end{pmatrix},
\]
where $z$ is and eigenvalue of $(2\beta I+ K^{-T}MK^{-1}M)z =\lambda z$.
\begin{proof}
Let $(\lambda,\bf{w})$ be an eigenpair of the matrix $\mathcal{P}^{-1}\mathcal{A}$, where $0\neq{\bf w}=(x;y;z)\in\Bbb{C}^{3m}$. Then, we have
${\mathcal{A}}{\bf w}=\lambda\mathcal{P} {\bf w}$, which is itself equivalent to
\begin{eqnarray}
2\beta  Mx - Mz&\h=\h& \lambda Ky, \label{EqEig1}\\
(1-\lambda) M y &\h=\h& -(1-\lambda) K^Tz,\label{EqEig2}\\
 (1-\lambda)Mx  &\h=\h& (1-\lambda) Ky \label{EqEig3}
 \end{eqnarray}
If $\lambda=1$, then obviously Eqs. (\ref{EqEig2}) and (\ref{EqEig3}) are always satisfied and Eq. (\ref{EqEig1}) takes the form
\[
x=\frac{1}{2\beta}(z+M^{-1}Ky) .
\]
Hence, every vector of the form \eqref{Eigw1} with $z$ or $y$ being nonzero is an eigenvector of $\mathcal{P}^{-1}{\mathcal{A}}$ corresponding to the eigenvalue $\lambda=1$.
Obviously, the algebraic multiplicity of $\lambda=1$ is equal to $2m$. We now assume that $\lambda\neq 1$. In this case, Eqs.  \eqref{EqEig1}-\eqref{EqEig3} can be rewritten in the form
\begin{eqnarray}
y &\h=\h& -M^{-1}K^Tz, \label{EqEigP1}\\
x &\h=\h& M^{-1}Ky = -M^{-1} KM^{-1}K^Tz,\label{EqEigP2}\\
 \lambda Ky&\h=\h&2\beta Mx-Mz . \label{EqEigP3}
\end{eqnarray}
Obviously, $z\neq 0$. Otherwise, from  Eqs.  \eqref{EqEigP1} and \eqref{EqEigP2} we have $x=y=z=0$ which is a contradiction with ${\bf w}=(x;y;z)$ being an eigenvector.  Substituting $x$ and $y$ from Eqs. \eqref{EqEigP1} and \eqref{EqEigP2} in Eq. \eqref{EqEigP3} yields
\[
-2\beta KM^{-1}K^Tz-Mz =-\lambda KM^{-1}K^Tz
\]
Premultiplying both sides of this equation by $z^*$ gives Eq. (\ref{Eignonunit}) and from Eqs. \eqref{EqEigP1} and \eqref{EqEigP2} the corresponding eigenvector is obtained.  Finally, from the above equation we have $(2\beta I+ K^{-T}MK^{-1}M)z =\lambda z$.	
\end{proof}
\end{thm}

In our numerical tests, we use the bilinear quadrilateral $\bf{Q}_1$ finite elements to discretize Eqs. \eqref{Eq001}-\eqref{Eq0003}.
Let $\mathcal{T}_h=\{\square_k:k=1,\ldots,K\}$ be a shape regular quadrilateral elements  of the domain $\Omega$. A sequence of quadrilateral grids  is said to be quasi-uniform if there exists a constant $\rho >0$ such that $\underline{h}\geq \rho h $,  where $\underline{h}=\min_{\square_k\in \tau_h}h_k$ and $h=\max_{\square_k\in \tau_h}h_k$ in which $h_k$ is denotes the diameter of $\square_k$, i.e., $h_k=\max\{|x-y|: x,y\in\overline{\square_k} \}$. We use the following theorem to obtain more concrete bound which are dependent on the PDE of the problem and the used finite element.
\begin{thm}\label{BoundMK} (\cite{Elman})
For $\bf{Q}_1$  approximation on a quasi-uniform subdivision of $\Bbb{R}^2$ for which a shape regularity condition holds, the mass matrix $M$ approximates the scaled identity matrix in the snese that
\[
c_1h^2 \leq \dfrac{\bx^TM\bx}{\bx^T\bx} \leq c_2 h^2,
\]
for every nonzero vector  ${\bf x}\in\Bbb{R}^n$. The constants $c_1$ and $c_2$ are independent of $h$ and $\beta$. On the other hand, the Galerkin matrix $K$ satisfies
\[
d_1h^2 \leq \dfrac{\bx^TK\bx}{\bx^T\bx} \leq d_2,
\]
for every nonzero vector  ${\bf x}\in\Bbb{R}^n$. The constants $d_1$ and $d_2$ are positive and independent of $h$ and $\beta$.
\end{thm}

Here, for $\bf{Q}_1$ approximation,  we mention that the mass matrix  $M$ and the stiffness matrix $K$ are symmetric positive definite.  From Theorem \ref{BoundMK}, it is easy to see that
\begin{eqnarray}
\frac{1 }{c_2h^2}&\h\leq\h&\frac{x^TM^{-1}x}{x^Tx}\leq \frac{1}{c_1h^2},\label{XMIX}\\
d_1h^2&\h\leq\h&\frac{x^Tx}{x^TK^{-1}x}\leq d_2,\label{XKIX}
\end{eqnarray}
for every nonzero $x\in\Bbb{R}^m$. Now, the following theorem can be stated.

\begin{thm}\label{EigBound}
For the $\bf{Q}_1$ approximation the non-unit eigenvalues of the matrix $\mathcal{P}^{-1}{\mathcal{A}}$ satisfies
\begin{equation}\label{NonUnitBound}
2\beta+\frac{c_1^2 h^4}{d_2^2} \leq \lambda\leq 2\beta+\frac{c_2^2}{d_1^2  },
\end{equation}
where $c_1$, $c_2$, $d_1$ and $d_2$ were defined in Theorem \ref{BoundMK}.
\begin{proof}
Let $0 \neq x\in\Bbb{R}^{n}$. Then we have
\begin{eqnarray}
 \frac{x^TKM^{-1}K^Tx}{x^Tx}&\h=\h& \frac{y^TM^{-1}y}{y^Ty}\frac{x^TK^2x}{x^Tx}\qquad (\text{where } y=Kx) \nonumber\\
                            &\h=\h& \frac{y^TM^{-1}y}{y^Ty}\frac{(K^{\frac{1}{2}}x)^T K (K^{\frac{1}{2}} x)}{x^Tx} \nonumber\\
                            &\h=\h& \frac{y^TM^{-1}y}{y^Ty}\frac{w^T K w}{w^Tw} \frac{w^Tw}{w^T K^{-1} w} \qquad (\text{where } w=K^{\frac{1}{2}}x). \label{EqKMIM}
\end{eqnarray}
Now using Theorem \ref{BoundMK}, Eqs. \eqref{XMIX} and \eqref{XKIX} we get
\[
\frac{d_1^2h^2}{c_2}  \leq \frac{x^TKM^{-1}K^Tx}{x^Tx} \leq \frac{d_2^2}{c_1h^2}.
\]
Finally, from the latter equation along with \eqref{EqKMIM} the desired result is obtained.
\end{proof}
\end{thm}

From  Theorem  \ref{EigBound} we see that the non-unit eigenvalues of the preconditioned matrix are real. The upper bound  for the  non-unit eigenvalues of the  preconditioned matrix is independent of the mesh size, and since the mesh size is usually small, the lower  bound of the  non-unit eigenvalues of the  preconditioned matrix can be estimated by $2\beta$. This means that the interval $(2\beta,2\beta+ (c_2/d_1)^2 ]$  is a good approximation for the interval containing the non-unit eigenvalues of the preconditioned matrix, which is independent of the mesh size. According to these comments it is expected that convergence  of a Krylov method like GMRES would be independent of  the mesh size.

In two dimensions, for $\Omega=[0,1]\times [0,1]$ and  for $\bf{Q}_1$  square elements (as our test example), the Fourier analysis gives that $c_1=1/9$, $c_2=1$, $d_1=2 \pi^2$ and $d_2=4$ (see \cite{Rees2010}). In this case, substituting these values in  Eq. \eqref{NonUnitBound} gives the following bounds for the non-unit eigenvalues of the preconditioned matrix $\mathcal{P}^{-1}{\mathcal{A}}$
\begin{equation}\label{BoundPIA}
2\beta+\frac{1}{1296}h^4 \leq \lambda \leq  2\beta+\frac{1}{4\pi^4}.
\end{equation}
From Theorems \ref{Thm1}, \ref{EigBound} and the latter equation we observe that if
\[
\beta\leq \frac{1}{2} \left( 1-\frac{1}{4\pi^4}\right)\approx 0.4987,
\]
then the eigenvalues of the preconditioned matrix  $\mathcal{P}^{-1}{\mathcal{A}}$ are contained in the interval $(2\beta,1]$. Therefore, the Chebyshev acceleration method can be applied to the preconditioned system ${P}^{-1}\mathcal{A}{\bf x}=\mathcal{P}^{-1}{\bf g}$  (see \cite{Saadbook}). Our numerical results show that this method can not compete with the GMRES method for solving the preconditioned system. However, the Chebyshev acceleration method is quite suitable  when we solve the system in a parallel environment.

Implementation of the preconditioner $\mathcal{P}$ in a Krylov subspace method can be done as follows.
Let $r=(r_1;r_2;r_3)$ and $t=(x;y;z)=\mathcal{P}^{-1}r$. Then it is easy to see that  the vector $t$ can be computed via the following steps.
\begin{enumerate}
	\item Solve $Mx=r_1-r_3$ for $x$.  \\ [-.75cm]
	\item Solve $Ky=r_1$ for $y$ \\[-.75cm]
	\item Solve $K^Tz=r_2-My$ for $z$
\end{enumerate}

As we see, in the implementation of the proposed preconditioner three systems with the matrices $K$, $K^T$ and  $M$ should be solved. For  the $\bf{Q}_1$ approximation,  both $M$ and $K$ are symmetric positive definite and  the systems with the coefficient matrices $M$ and $K$ can be solved exactly using the Cholesky factorization or inexactly using the CG iterative method.

\section{Numerical experiments}\label{Sec3}

We consider the problem
\begin{eqnarray*}
	\min_{u,f}~~~ &\h\h&\frac{1}{2} \|u-u_*\|_{L^2(\Omega)}^2+\beta\|f\|_{L^2(\Omega)}^2,\label{Eq0001}\\
	\text{s.t.}~~~~&\h\h& -\Delta u=f\quad \text{in}~~ \Omega,\label{Eq0002}\\
	&\h\h&u=\hat{u}|_{\partial \Omega} \quad \text{on}~~\partial\Omega,\label{Eq003}
\end{eqnarray*}
where $\Omega=[0,1]\times [0,1]$  and
\begin{equation*}\label{Eq04}
	\hat{u}=\left\{
	\begin{array}{ll}
		(2x-1)^2(2y-1)^2 & \text{if } (x,y)\in [0,\frac{1}{2}]^2,   \\[2mm]
		0, & \text{elsewhere}.
	\end{array}
	\right.
\end{equation*}
We use the codes$^{1}$\footnote{https://github.com/tyronerees/poisson-control} of  the  paper \cite{Rees2010} to generate the system \eqref{MainSys}. The problem is discretized using bilinear quadrilateral $\bf{Q}_1$ and a grid of $2^{\ell} \times 2^{\ell}$ is used. In this case, $\mathcal{A}$ is an $n\times n $  matrix with $n=3(2^{\ell}-1)(2^{\ell}-1)$.

\begin{figure}[!ht]
	\centering
	\subfigure{\label{eigAh}
		\includegraphics*[width=.45\textwidth]{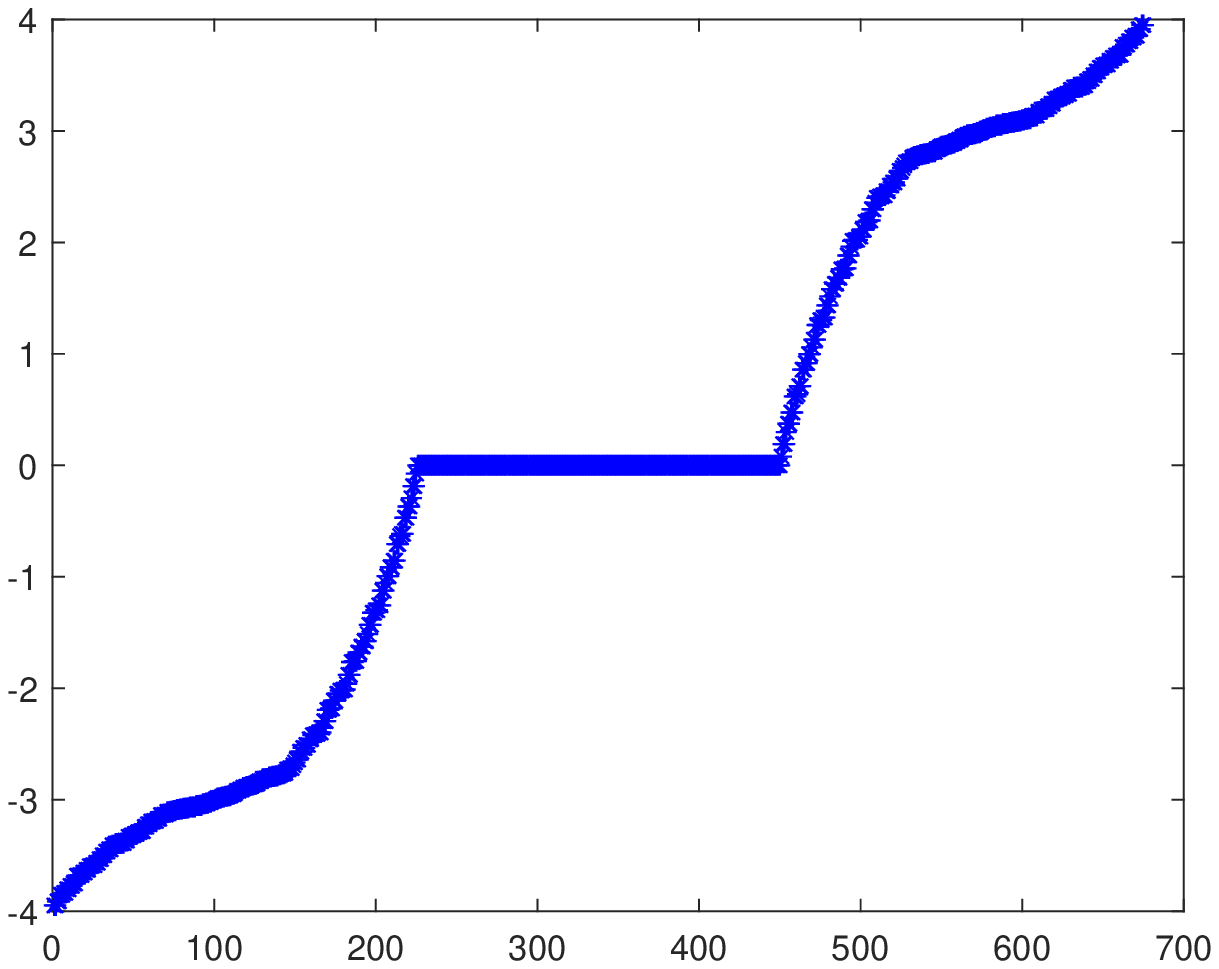}}
	\hspace{3mm}
	\subfigure{\label{mylabel1}
		\includegraphics*[width=.45\textwidth]{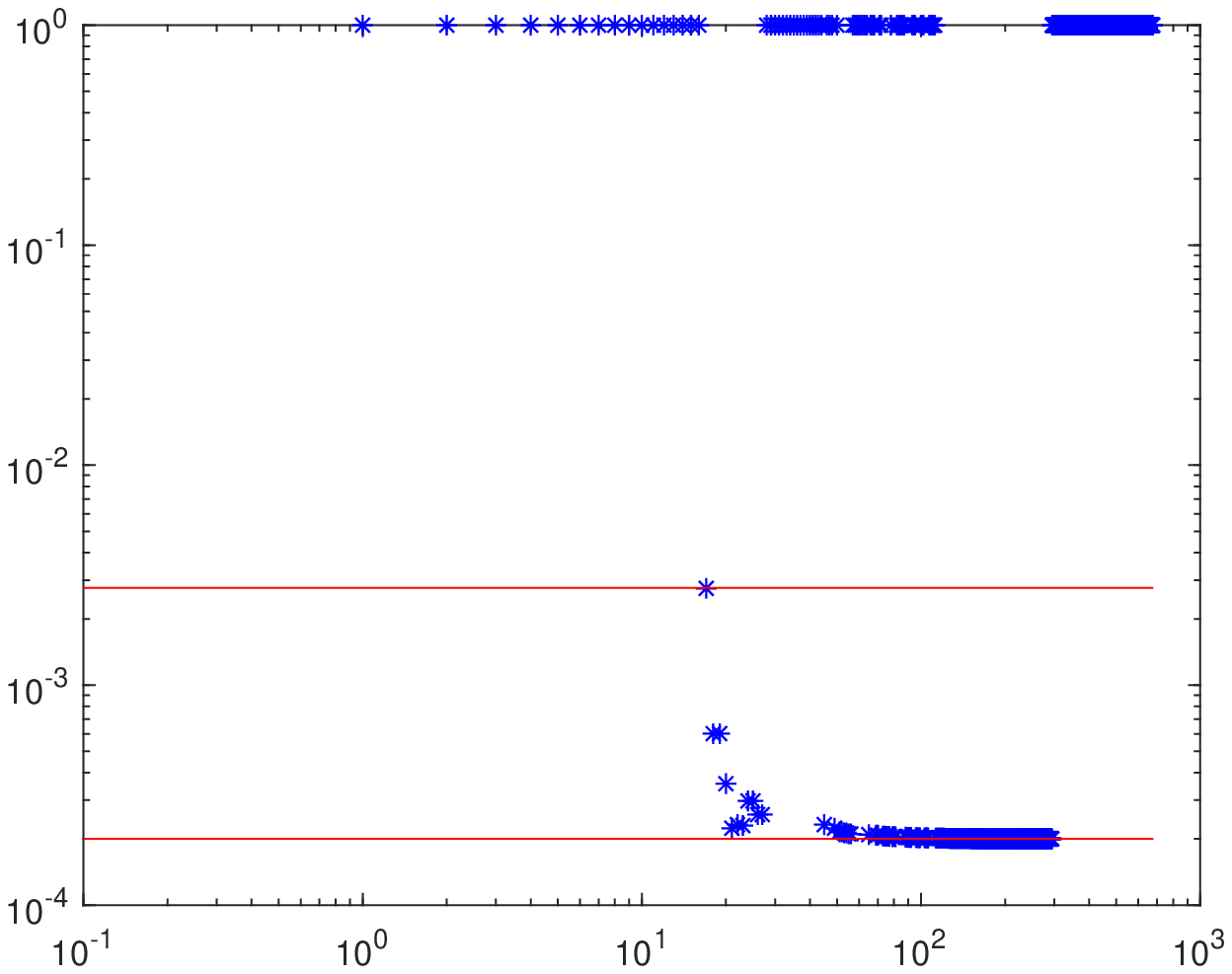}}
	\caption{Eigenvalues distribution of the original matrix $\mathcal{A}$ (left) and the preconditioned matrix $\mathcal{P}^{-1}\mathcal{A}$ (right) for $(h,\beta)=(2^{-4},10^{-4})$. This right figure has been displayed in log-log format.}
	\label{Fig1}
\end{figure}

We first display the eigenvalue distribution of the matrices  $\mathcal{A}$ and  $\mathcal{P}^{-1}\mathcal{A}$ for $(h,\beta)=(2^{-4},10^{-4})$. To do so, in Figure \ref{Fig1}  we display ${\{(i,\lambda_i)\}}_{i=1}^n$ for these matrices, where $\lambda_i$, $i=1,\ldots,n,$ are the eigenvalues of the matrices.  
The lower and upper bounds of the non-unit eigenvalues of $\mathcal{P}^{-1}\mathcal{A}$ given in Eq. \eqref{BoundPIA} have been displayed by two horizontal  solid  lines.   As we see, Figure \ref{Fig1} confirms  Eq. \ref{BoundMK}. 

In the sequel, we present two sets of the numerical  results. In the first set, we use the complete version of the preconditioned GMRES method for solving \eqref{MainSys} and in this case all the sub-systems are solved exactly using the Cholesky factorization (GMRES/Cholesky). We use right-preconditioning.
Numerical results are compared in terms of  both the number of iterations and the CPU time which are, respectively, denoted by ``IT" and ``CPU" in the tables. In all the tests, we use a zero vector as an initial guess and the iteration  is  stopped as soon as the residual norm  is reduced by a factor of $10^{6}$. The maximum number of the iterations is set to be $maxit=\min\{500,3m\}$. A ``--" means that  the method has not converged in $maxit$ iterations.
All runs are implemented in \textsc{Matlab} R2015 , equipped with a Laptop with 2.60 GHz central processing unit (Intel(R) Core(TM) i7-4510), 8 GB memory and Windows 8 operating system.

We compare the numerical results of the proposed preconditioner (denoted by $\mathcal{P}$ in the tables) with those of the preconditioners reviewed in Section \ref{Sec2}. The first set of the numerical results have been presented in Table \ref{Tb1}. 
We have highlighted the best result for each choice of the parameter in boldface.
As we see the GMRES method in conjunction with the proposed preconditioner converges in a small  number of the iterations for different values of the step size and the regularization parameter. This is not the case for the other methods. For example the number of iterations for the preconditioner $\mathcal{P}_D$ increases when the regularization parameter becomes small. The same observations can be made also for the other methods.

For the second set of the numerical experiments, we  apply the flexible GMRES (FGMRES) method (see \cite{FGMRES}) for solving the preconditioned system and in this case  the sub-systems (for all the preconditioners) are solved inexactly using the preconditioned CG method (FGMRES/PCG).  We use the incomplete Cholesky factorization with dropping tolerance $\tau=0.01$ as the preconditioner for the sub-systems. To do so, we have used the following command of \textsc{Matlab}:
\begin{center}
\verb"L = ichol(M,struct('type','ict','droptol',1e-02))".	
\end{center}
Left-preconditioning technique is used for the subsystems.  The inner iteration is terminated as soon as the residual norm is reduced by a  factor of $10^{3}$.  The  maximum number of iterations of the PCG method is set to be $\min\{m,20\}$. All the other assumptions are similar to those of the first set of the numerical results. Numerical results have been given in Table \ref{Tb2}. This table shows the superiority of the proposed preconditioner to the other ones. In many cases the other preconditioners fail to converge, whereas the proposed preconditioner never fail for the tested systems. As the numerical results show there is no significant difference between the two set of the numerical results.
	\begin{table}
		\centering\caption{Numerical results for GMRES/Cholesky.\label{Tb1}}
		\begin{center}
\scriptsize
\scalebox{.7}
{
\begin{tabular}{|c|c|c|c|c|c|c|c|c|c|c|c|}
\hline
&&\multicolumn{1}{c|}{$\mathcal{P}$}&\multicolumn{1}{c|}{D}&\multicolumn{1}{c|}{BCD}&\multicolumn{1}{c|}{BCT}&\multicolumn{1}{c|}{C}&\multicolumn{1}{c|}{BS}&\multicolumn{1}{c|}{BLT}&\multicolumn{1}{c|}{ID}&\multicolumn{1}{c|}{P1/P2}&\multicolumn{1}{c|}{P3/P4}\\ \hline
					$\beta$& h &IT(CPU) &IT(CPU) &IT(CPU) &IT(CPU) &IT(CPU) &IT(CPU) &IT(CPU) &IT(CPU) &IT(CPU)/IT(CPU) &IT(CPU)/IT(CPU)\\ \hline
$10^{-1}$&$2^{-2}$&4(0.01)&5(0.01)&18(0.01)&7(0.01)&8(0.01)&12(0.01)&7(0.01)&3(0.01)&3(0.00)/\bf2(0.01)&7(0.01)/7(0.01) \\
&$2^{-3}$&\bf3(0.00)&6(0.00)&98(0.06)&33(0.01)&20(0.01)&64(0.03)&34(0.01)&4(0.00)&\bf3(0.00)/3(0.00)&33(0.01)/36(0.01) \\
&$2^{-4}$&3(0.00)&5(0.00)&-(-)&165(0.21)&27(0.02)&330(0.75)&171(0.24)&3(0.00)&3(0.00)/\bf2(0.00)&171(0.22)/181(0.26) \\
&$2^{-5}$&3(0.02)&5(0.02)&-(-)&-(-)&31(0.10)&-(-)&-(-)&3(0.02)&3(0.01)/\bf2(0.01)&-(-)/-(-) \\
&$2^{-6}$&3(0.13)&5(0.18)&-(-)&-(-)&33(0.87)&-(-)&-(-)&3(0.15)&3(0.12)/\bf2(0.11)&-(-)/-(-) \\
&$2^{-7}$&3(1.08)&4(1.33)&-(-)&-(-)&35(7.89)&-(-)&-(-)&3(1.17)&2(0.83)/\bf2(0.83)&-(-)/-(-) \\
\hline
$10^{-2}$&$2^{-2}$&4(0.01)&7(0.01)&18(0.01)&7(0.01)&8(0.01)&12(0.01)&7(0.01)&4(0.01)&\bf4(0.01)/4(0.01)&7(0.01)/7(0.01) \\
&$2^{-3}$&4(0.00)&7(0.00)&97(0.06)&32(0.01)&18(0.01)&64(0.03)&33(0.01)&4(0.00)&4(0.00)/\bf3(0.00)&33(0.01)/36(0.01) \\
&$2^{-4}$&4(0.00)&7(0.01)&-(-)&163(0.20)&24(0.01)&334(0.74)&168(0.21)&4(0.00)&4(0.00)/\bf3(0.00)&169(0.21)/181(0.24) \\
&$2^{-5}$&4(0.02)&7(0.02)&-(-)&-(-)&28(0.08)&-(-)&-(-)&4(0.02)&\bf3(0.01)/3(0.01)&-(-)/-(-) \\
&$2^{-6}$&\bf3(0.12)&7(0.21)&-(-)&-(-)&30(0.77)&-(-)&-(-)&4(0.24)&\bf3(0.16)/3(0.15)&-(-)/-(-) \\
&$2^{-7}$&\bf3(1.32)&6(1.76)&-(-)&-(-)&32(7.34)&-(-)&-(-)&4(1.35)&\bf3(0.97)/3(0.97)&-(-)/-(-) \\
\hline
$10^{-3}$&$2^{-2}$\bf&5(0.01)&9(0.01)&18(0.01)&7(0.01)&8(0.01)&12(0.01)&7(0.01)&6(0.01)&\bf5(0.01)/5(0.01)&7(0.01)/7(0.01) \\
&$2^{-3}$&6(0.00)&11(0.00)&97(0.06)&31(0.01)&15(0.01)&64(0.03)&33(0.01)&6(0.00)&\bf5(0.00)/5(0.00)&33(0.01)/35(0.01) \\
&$2^{-4}$&6(0.00)&10(0.01)&-(-)&160(0.19)&21(0.01)&336(0.75)&165(0.21)&6(0.01)&\bf5(0.00)/5(0.00)&165(0.20)/180(0.24) \\
&$2^{-5}$&\bf4(0.02)&9(0.03)&-(-)&-(-)&25(0.07)&-(-)&-(-)&5(0.02)&5(0.02)/5(0.02)&-(-)/-(-) \\
&$2^{-6}$&\bf4(0.14)&9(0.24)&-(-)&-(-)&28(0.73)&-(-)&-(-)&5(0.17)&\bf4(0.13)/4(0.13)&-(-)/-(-) \\
&$2^{-7}$&\bf3(1.00)&8(1.91)&-(-)&-(-)&30(6.69)&-(-)&-(-)&5(1.48)&4(1.07)/4(1.13)&-(-)/-(-) \\
\hline
$10^{-4}$&$2^{-2}$&\bf6(0.02)&13(0.02)&18(0.02)&7(0.02)&8(0.03)&12(0.02)&7(0.02)&7(0.03)&7(0.02)/7(0.03)&7(0.02)/7(0.01) \\
&$2^{-3}$&\bf7(0.01)&15(0.01)&95(0.11)&27(0.01)&16(0.02)&64(0.06)&33(0.02)&9(0.00)&$\bf7(0.00)$/8(0.00)&32(0.02)/34(0.02) \\
&$2^{-4}$&\bf7(0.02)&15(0.02)&448(1.75)&97(0.09)&20(0.01)&331(0.83)&163(0.22)&8(0.01)&\bf7(0.01)/7(0.01)&156(0.21)/174(0.27) \\
&$2^{-5}$&\bf6(0.02)&15(0.04)&-(-)&287(2.89)&23(0.08)&-(-)&-(-)&8(0.03)&7(0.02)/7(0.02)&-(-)/-(-) \\
&$2^{-6}$&\bf6(0.25)&15(0.47)&-(-)&-(-)&27(0.85)&-(-)&-(-)&8(0.28)&7(0.21)/7(0.21)&-(-)/-(-) \\
&$2^{-7}$&\bf4(1.42)&13(3.21)&-(-)&153(26.79)&29(7.80)&-(-)&-(-)&7(2.26)&6(1.72)/6(1.69)&-(-)/-(-) \\
\hline
$10^{-5}$&$2^{-2}$&8(0.01)&13(0.02)&15(0.01)&\bf4(0.01)&8(0.02)&12(0.01)&7(0.01)&7(0.02)&7(0.02)/7(0.02)&6(0.02)/7(0.01) \\
&$2^{-3}$&10(0.01)&25(0.02)&58(0.05)&\bf9(0.01)&17(0.01)&46(0.03)&27(0.02)&15(0.01)&11(0.01)/14(0.01)&19(0.01)/27(0.01) \\
&$2^{-4}$&\bf10(0.01)&27(0.02)&248(0.52)&23(0.01)&23(0.02)&194(0.31)&115(0.13)&15(0.01)&11(0.01)/14(0.01)&75(0.06)/128(0.15) \\
&$2^{-5}$&\bf8(0.02)&25(0.07)&-(-)&28(0.07)&26(0.08)&-(-)&435(6.22)&13(0.05)&11(0.03)/12(0.03)&262(2.25)/-(-) \\
&$2^{-6}$&\bf7(0.23)&23(0.65)&-(-)&26(0.52)&28(0.93)&-(-)&-(-)&12(0.36)&10(0.27)/11(0.27)&-(-)/-(-) \\
&$2^{-7}$&\bf6(1.76)&21(4.86)&-(-)&24(3.99)&31(8.75)&-(-)&-(-)&11(3.07)&10(2.19)/11(2.32)&-(-)/-(-) \\
\hline
$10^{-6}$&$2^{-2}$&8(0.01)&13(0.02)&8(0.01)&\bf2(0.01)&8(0.01)&8(0.01)&5(0.01)&7(0.01)&7(0.01)/7(0.02)&3(0.01)/5(0.01) \\
&$2^{-3}$&12(0.01)&43(0.03)&23(0.01)&\bf2(0.00)&22(0.01)&21(0.01)&13(0.01)&24(0.01)&17(0.01)/24(0.01)&8(0.00)/14(0.01) \\
&$2^{-4}$&12(0.01)&50(0.05)&86(0.08)&\bf3(0.00)&29(0.02)&80(0.07)&48(0.04)&29(0.02)&21(0.02)/26(0.02)&27(0.02)/53(0.04) \\
&$2^{-5}$&11(0.03)&49(0.16)&344(4.06)&\bf3(0.02)&35(0.12)&313(3.20)&184(1.40)&26(0.09)&20(0.06)/25(0.07)&90(0.37)/220(1.65) \\
&$2^{-6}$&10(0.32)&47(1.28)&-(-)&\bf2(0.13)&39(1.25)&-(-)&-(-)&23(0.62)&19(0.45)/22(0.54)&264(9.50)/-(-) \\
&$2^{-7}$&10(2.63)&43(9.69)&-(-)&\bf2(0.88)&42(11.69)&-(-)&-(-)&22(5.17)&18(3.60)/21(3.84)&-(-)/-(-) \\
\hline
$10^{-7}$&$2^{-2}$&8(0.02)&13(0.02)&5(0.01)&\bf1(0.01)&7(0.02)&5(0.01)&4(0.01)&7(0.02)&7(0.02)/7(0.02)&2(0.02)/4(0.01) \\
&$2^{-3}$&12(0.01)&56(0.04)&8(0.00)&\bf1(0.00)&27(0.01)&10(0.01)&7(0.01)&31(0.02)&25(0.01)/33(0.02)&3(0.00)/7(0.01) \\
&$2^{-4}$&12(0.01)&90(0.11)&26(0.01)&\bf1(0.00)&47(0.03)&32(0.02)&19(0.01)&53(0.04)&34(0.02)/49(0.04)&9(0.01)/20(0.01) \\
&$2^{-5}$&10(0.03)&96(0.42)&101(0.43)&\bf1(0.02)&58(0.22)&115(0.52)&70(0.26)&56(0.25)&37(0.13)/51(0.19)&31(0.10)/80(0.33) \\
&$2^{-6}$&5(0.19)&89(2.64)&386(16.34)&\bf1(0.09)&66(2.03)&403(17.71)&251(8.40)&53(1.41)&36(0.84)/48(1.06)&90(1.98)/317(12.10) \\
&$2^{-7}$&2(1.15)&82(18.46)&-(-)&\bf1(0.70)&72(18.60)&-(-)&-(-)&48(10.11)&34(5.67)/43(7.07)&213(39.67)/-(-) \\
\hline
$10^{-8}$&$2^{-2}$&8(0.01)&15(0.02)&5(0.01)&\bf1(0.01)&7(0.01)&4(0.01)&3(0.01)&7(0.02)&7(0.01)/8(0.02)&2(0.01)/3(0.01) \\
&$2^{-3}$&12(0.01)&58(0.04)&5(0.00)&\bf1(0.00)&26(0.01)&6(0.00)&4(0.00)&32(0.02)&30(0.01)/35(0.02)&2(0.00)/4(0.00) \\
&$2^{-4}$&12(0.02)&136(0.23)&8(0.01)&\bf1(0.01)&49(0.04)&10(0.01)&8(0.01)&78(0.07)&52(0.04)/78(0.07)&3(0.01)/8(0.01) \\
&$2^{-5}$&8(0.03)&175(1.13)&26(0.06)&\bf1(0.01)&94(0.46)&39(0.14)&25(0.08)&101(0.56)&64(0.26)/96(0.48)&9(0.03)/27(0.08) \\
&$2^{-6}$&5(0.21)&176(6.50)&107(2.44)&\bf1(0.09)&115(3.92)&146(3.74)&90(2.02)&102(2.95)&65(1.40)/96(2.33)&30(0.60)/106(2.50) \\
&$2^{-7}$&2(1.19)&164(37.22)&398(90.13)&\bf1(0.71)&130(34.23)&-(-)&312(64.68)&96(20.48)&63(10.39)/89(15.10)&78(12.28)/411(93.96) \\
\hline
$10^{-9}$&$2^{-2}$&8(0.01)&15(0.02)&2(0.01)&\bf1(0.01)&7(0.02)&4(0.01)&3(0.01)&7(0.01)&7(0.01)/9(0.02)&$\bf1(0.01)$/3(0.01) \\
&$2^{-3}$&12(0.01)&60(0.04)&2(0.00)&\bf1(0.00)&27(0.01)&4(0.00)&3(0.00)&32(0.02)&33(0.01)/39(0.02)&2(0.00)/3(0.00) \\
&$2^{-4}$&12(0.01)&154(0.23)&5(0.01)&\bf1(0.00)&58(0.05)&6(0.01)&5(0.01)&92(0.09)&65(0.05)/100(0.10)&2(0.01)/5(0.01) \\
&$2^{-5}$&8(0.03)&284(2.67)&5(0.02)&\bf1(0.01)&133(0.77)&13(0.03)&9(0.03)&167(1.05)&96(0.41)/163(0.99)&3(0.02)/10(0.03) \\
&$2^{-6}$&5(0.19)&343(15.36)&23(0.43)&\bf1(0.09)&188(7.31)&48(0.96)&31(0.64)&199(7.07)&113(2.99)/184(5.54)&9(0.23)/34(0.70) \\
&$2^{-7}$&2(1.11)&354(96.31)&122(20.11)&\bf1(0.70)&222(65.87)&189(34.13)&111(18.34)&210(50.55)&111(19.16)/182(34.05)&26(4.25)/137(23.45) \\
\hline
$10^{-10}$&$2^{-2}$&8(0.01)&15(0.02)&2(0.01)&\bf1(0.01)&7(0.01)&2(0.01)&3(0.01)&7(0.02)&7(0.01)/10(0.02)&$\bf1(0.01)$/3(0.01) \\
&$2^{-3}$&12(0.01)&60(0.04)&2(0.00)&\bf1(0.00)&29(0.01)&4(0.00)&3(0.00)&32(0.01)&33(0.01)/41(0.02)&$\bf1(0.00)$/3(0.00) \\
&$2^{-4}$&12(0.01)&158(0.22)&2(0.00)&\bf1(0.00)&83(0.08)&4(0.01)&4(0.00)&98(0.10)&74(0.06/112(0.14)&$\bf1(0.01)$/4(0.00) \\
&$2^{-5}$&8(0.03)&328(3.88)&2(0.01)&\bf1(0.01)&159(1.08)&6(0.03)&5(0.02)&216(1.96)&125(0.67)/228(2.00)&2(0.02)/5(0.02) \\
&$2^{-6}$&5(0.20)&-(-)&5(0.16)&\bf1(0.09)&271(12.57)&16(0.32)&15(0.39)&364(17.65)&170(4.83)/331(13.78)&2(0.12)/11(0.24) \\
&$2^{-7}$&2(1.08)&-(-)&41(6.77)&\bf1(0.68)&432(147.07)&63(9.76)&38(5.87)&494(154.54)&187(35.72)/377(91.84)&7(1.60)/41(7.97) \\
\hline
\end{tabular}
}
\end{center}
\end{table}

\begin{table}
		\centering\caption{Numerical results for FGMRES/PCG. \label{Tb2}}
		\begin{center}
			\label{ibexact}
			\scriptsize
			\scalebox{.75}
			{
				\begin{tabular}{|c|c|c|c|c|c|c|c|c|c|c|c|}
					\hline
					&&\multicolumn{1}{c|}{$\mathcal{P}$}&\multicolumn{1}{c|}{D}&\multicolumn{1}{c|}{BCD}&\multicolumn{1}{c|}{BCT}&\multicolumn{1}{c|}{C}&\multicolumn{1}{c|}{BS}&\multicolumn{1}{c|}{BLT}&\multicolumn{1}{c|}{ID}&\multicolumn{1}{c|}{P1/P2}&\multicolumn{1}{c|}{P3/P4}\\ \hline
					$\beta$& h &IT(CPU) &IT(CPU) &IT(CPU) &IT(CPU) &IT(CPU) &IT(CPU) &IT(CPU) &IT(CPU) &IT(CPU)/IT(CPU) &IT(CPU)/IT(CPU)\\ \hline
					$10^{-1}$&$2^{-2}$&4(0.01)&5(0.04)&18(0.04)&7(0.02)&8(0.05)&13(0.03)&7(0.02)&3(0.01)& 3(0.01)/$\bf2(0.01)$&7(0.02)/7(0.02) \\
&$2^{-3}$&5(0.01)&7(0.02)&96(0.21)&33(0.06)&35(0.11)&65(0.13)&34(0.06)&5(0.01)&$\bf3(0.01)/3(0.01)$&34(0.06)/35(0.06) \\
&$2^{-4}$&4(0.01)&7(0.02)&-(-)&164(0.53)&236(1.24)&339(1.68)&171(0.57)&5(0.01)&$\bf3(0.01)/3(0.01)$&170(0.56)/178(0.60) \\
&$2^{-5}$&4(0.02)&7(0.03)&-(-)&-(-)&-(-)&-(-)&-(-)&5(0.02)&$\bf3(0.01)/3(0.01)$&-(-)/-(-) \\
&$2^{-6}$&4(0.07)&7(0.12)&-(-)&-(-)&-(-)&-(-)&-(-)&5(0.11)&3(0.06)/$\bf2(0.04)$&-(-)/-(-) \\
&$2^{-7}$&6(0.45)&13(1.17)&-(-)&-(-)&-(-)&-(-)&-(-)&13(1.04)&$\bf3(0.24)/3(0.24)$&-(-)/-(-) \\
\hline
$10^{-2}$&$2^{-2}$&$\bf4(0.01)$&7(0.01)&18(0.03)&7(0.01)&8(0.02)&13(0.02)&7(0.01)&$\bf4(0.01)$&$\bf4(0.01)/4(0.01)$&7(0.01)/7(0.01) \\
&$2^{-3}$&6(0.01)&9(0.02)&96(0.21)&32(0.06)&29(0.09)&65(0.13)&33(0.06)&6(0.01)&4(0.01)/$\bf3(0.01)$&33(0.06)/34(0.06) \\
&$2^{-4}$&6(0.01)&7(0.02)&-(-)&162(0.52)&122(0.52)&339(1.68)&168(0.55)&5(0.01)&4(0.01)/$\bf3(0.01)$&168(0.55)/175(0.58) \\
&$2^{-5}$&6(0.03)&8(0.04)&-(-)&-(-)&-(-)&-(-)&-(-)&6(0.03)&4(0.02)/$\bf3(0.01)$&-(-)/-(-) \\
&$2^{-6}$&4(0.07)&9(0.17)&-(-)&-(-)&-(-)&-(-)&-(-)&8(0.14)&$\bf3(0.05)/3(0.05)$&-(-)/-(-) \\
&$2^{-7}$&6(0.55)&16(1.50)&-(-)&-(-)&-(-)&-(-)&-(-)&15(1.22)&$\bf3(0.24)/3(0.25)$&-(-)/-(-) \\
\hline
$10^{-3}$&$2^{-2}$&$\bf5(0.01)$&9(0.02)&18(0.03)&7(0.01)&8(0.02)&13(0.02)&7(0.01)&6(0.02)&$\bf5(0.01)/5(0.01)$&7(0.01)/7(0.01) \\
&$2^{-3}$&8(0.02)&11(0.03)&96(0.21)&31(0.06)&20(0.06)&65(0.13)&32(0.06)&8(0.02)&$\bf5(0.01)/5(0.01)$&33(0.06)/34(0.06) \\
&$2^{-4}$&8(0.02)&11(0.03)&-(-)&160(0.58)&94(0.36)&337(1.67)&165(0.54)&7(0.02)&$\bf5(0.01)/5(0.01)$&165(0.55)/171(0.58) \\
&$2^{-5}$&7(0.03)&11(0.06)&-(-)&-(-)&-(-)&-(-)&-(-)&8(0.04)&$\bf5(0.02)/5(0.02)$&-(-)/-(-) \\
&$2^{-6}$&7(0.14)&11(0.21)&-(-)&-(-)&-(-)&-(-)&-(-)&11(0.21)&$\bf4(0.08)/4(0.07)$&-(-)/-(-) \\
&$2^{-7}$&9(0.71)&18(1.73)&-(-)&-(-)&-(-)&-(-)&-(-)&55(5.30)&$\bf5(0.41)/5(0.42)$&-(-)/-(-) \\
\hline
$10^{-4}$&$2^{-2}$&$\bf6(0.01)$&13(0.03)&18(0.03)&7(0.01)&8(0.02)&13(0.02)&7(0.01)&7(0.01)&7(0.01)/7(0.01)&7(0.01)/7(0.01) \\
&$2^{-3}$&11(0.02)&17(0.04)&92(0.20)&27(0.05)&24(0.07)&64(0.13)&32(0.06)&12(0.03)&$\bf7(0.01)$/9(0.02)&31(0.05)/32(0.06) \\
&$2^{-4}$&11(0.03)&17(0.05)&445(2.67)&97(0.25)&66(0.23)&329(1.59)&161(0.51)&12(0.03)&$\bf8(0.02)/8(0.02)$&155(0.48)/165(0.53) \\
&$2^{-5}$&9(0.04)&17(0.09)&-(-)&287(2.76)&-(-)&-(-)&-(-)&13(0.07)&$\bf8(0.03)/8(0.03)$&-(-)/-(-) \\
&$2^{-6}$&9(0.16)&16(0.31)&-(-)&-(-)&-(-)&-(-)&-(-)&17(0.34)&$\bf7(0.13)$/8(0.14)&-(-)/-(-) \\
&$2^{-7}$&10(0.80)&28(2.57)&-(-)&153(9.50)&-(-)&-(-)&-(-)&76(8.21)&9(0.75)/$\bf8(0.66)$&-(-)/-(-) \\
\hline
$10^{-5}$&$2^{-2}$&8(0.01)&13(0.03)&15(0.02)&4(0.01)&8(0.02)&11(0.02)&7(0.01)&7(0.02)&7(0.01)/7(0.01)&$\bf6(0.01)$/7(0.01) \\
&$2^{-3}$&15(0.03)&28(0.07)&58(0.11)&9(0.01)&21(0.06)&43(0.08)&27(0.05)&21(0.05)&$\bf14(0.03)$/16(0.03)&19(0.03)/27(0.05) \\
&$2^{-4}$&16(0.04)&29(0.09)&247(1.00)&23(0.04)&31(0.10)&218(0.82)&115(0.32)&23(0.07)&$\bf14(0.03)$/15(0.04)&74(0.17)/126(0.36) \\
&$2^{-5}$&16(0.08)&29(0.15)&-(-)&28(0.08)&-(-)&-(-)&435(5.99)&24(0.13)&15(0.06)/$\bf14(0.06)$&262(2.38)/-(-) \\
&$2^{-6}$&$\bf12(0.22)$&27(0.51)&-(-)&26(0.19)&-(-)&-(-)&-(-)&28(0.60)&15(0.27)/13(0.23)&-(-)/-(-) \\
&$2^{-7}$&$\bf13(1.06)$&46(4.46)&-(-)&24(0.51)&-(-)&-(-)&-(-)&-(-)&43(4.00)/13(1.08)&-(-)/-(-) \\
\hline
$10^{-6}$&$2^{-2}$&8(0.01)&13(0.03)&8(0.01)&$\bf2(0.00)$&8(0.02)&7(0.01)&5(0.01)&7(0.01)&7(0.01)/7(0.01)&3(0.00)/5(0.01) \\
&$2^{-3}$&21(0.04)&46(0.12)&23(0.04)&$\bf2(0.00)$&25(0.07)&21(0.03)&13(0.02)&35(0.09)&29(0.06)/26(0.05)&8(0.01)/14(0.02) \\
&$2^{-4}$&22(0.06)&53(0.17)&86(0.21)&$\bf3(0.00)$&31(0.11)&88(0.22)&48(0.10)&54(0.18)&37(0.09)/29(0.07)&27(0.05)/52(0.11) \\
&$2^{-5}$&21(0.10)&53(0.29)&344(3.82)&$\bf3(0.01)$&46(0.26)&378(4.56)&184(1.29)&63(0.41)&46(0.23)/28(0.13)&90(0.40)/219(1.74) \\
&$2^{-6}$&19(0.37)&51(0.97)&-(-)&$\bf2(0.01)$&-(-)&-(-)&-(-)&472(27.78)&49(0.98)/25(0.43)&263(7.07)/-(-) \\
&$2^{-7}$&23(1.94)&70(7.11)&-(-)&$\bf2(0.03)$&-(-)&-(-)&-(-)&-(-)&202(28.64)/25(1.94)&-(-)/-(-) \\
\hline
$10^{-7}$&$2^{-2}$&8(0.01)&13(0.03)&5(0.01)&$\bf1(0.00)$&7(0.02)&4(0.01)&4(0.01)&7(0.01)&7(0.01)/7(0.01)&2(0.00)/4(0.01) \\
&$2^{-3}$&23(0.05)&58(0.15)&8(0.01)&$\bf1(0.00)$&31(0.09)&8(0.01)&7(0.01)&53(0.14)&63(0.14)/37(0.08)&3(0.00)/7(0.01) \\
&$2^{-4}$&22(0.05)&92(0.32)&26(0.05)&$\bf1(0.00)$&50(0.18)&33(0.06)&19(0.03)&214(1.02)&127(0.41)/52(0.14)&9(0.02)/20(0.04) \\
&$2^{-5}$&17(0.08)&98(0.64)&101(0.48)&$\bf1(0.00)$&62(0.38)&130(0.71)&69(0.28)&483(8.34)&179(1.50)/54(0.28)&31(0.09)/79(0.33) \\
&$2^{-6}$&8(0.15)&94(1.94)&386(14.33)&$\bf1(0.00)$&83(1.62)&-(-)&250(6.50)&-(-)&201(6.64)/50(0.94)&90(1.12)/315(9.85) \\
&$2^{-7}$&4(0.28)&109(10.73)&-(-)&$\bf1(0.02)$&104(8.85)&-(-)&-(-)&-(-)&403(83.66)/47(3.86)&213(17.97)/-(-) \\
\hline
$10^{-8}$&$2^{-2}$&8(0.01)&13(0.03)&5(0.01)&$\bf1(0.00)$&7(0.02)&4(0.01)&3(0.00)&7(0.01)&7(0.01)/7(0.01)&2(0.00)/3(0.00) \\
&$2^{-3}$&23(0.05)&59(0.16)&5(0.01)&$\bf1(0.00)$&29(0.08)&6(0.01)&4(0.01)&62(0.17)&96(0.23)/41(0.08)&2(0.00)/4(0.01) \\
&$2^{-4}$&23(0.06)&138(0.54)&8(0.01)&$\bf1(0.00)$&53(0.19)&10(0.02)&8(0.01)&331(1.95)&380(2.08)/84(0.25)&3(0.00)/8(0.01) \\
&$2^{-5}$&16(0.07)&177(1.50)&26(0.07)&$\bf1(0.00)$&98(0.67)&43(0.14)&25(0.07)&-(-)&-(-)/100(0.63)&9(0.02)/27(0.07) \\
&$2^{-6}$&7(0.13)&180(4.75)&107(1.48)&$\bf1(0.00)$&119(2.62)&168(3.18)&90(1.13)&-(-)&-(-)/99(2.23)&30(0.21)/106(1.46) \\
&$2^{-7}$&4(0.28)&194(21.83)&392(55.57)&$\bf1(0.01)$&162(16.24)&-(-)&312(35.98)&-(-)&-(-)/93(8.56)&78(3.04)/410(61.89) \\
\hline
$10^{-9}$&$2^{-2}$&8(0.01)&13(0.03)&2(0.00)&$\bf1(0.00)$&7(0.02)&4(0.01)&3(0.00)&7(0.01)&7(0.01)/7(0.01)&$\bf1(0.00)$/3(0.00) \\
&$2^{-3}$&23(0.05)&60(0.16)&2(0.00)&$\bf1(0.00)$&31(0.09)&4(0.01)&3(0.00)&74(0.21)&98(0.24)/44(0.10)&2(0.00)/3(0.00) \\
&$2^{-4}$&23(0.06)&154(0.61)&5(0.01)&$\bf1(0.00)$&60(0.22)&6(0.01)&5(0.01)&375(2.37)&436(2.63)/103(0.32)&2(0.00)/5(0.01) \\
&$2^{-5}$&16(0.08)&281(3.07)&5(0.01)&$\bf1(0.00)$&136(1.06)&12(0.03)&9(0.02)&-(-)&-(-)/165(1.39)&3(0.01)/10(0.02) \\
&$2^{-6}$&7(0.14)&329(13.20)&23(0.15)&$\bf1(0.00)$&190(5.22)&54(0.52)&31(0.23)&-(-)&-(-)/185(5.66)&9(0.05)/34(0.27) \\
&$2^{-7}$&4(0.29)&362(59.55)&104(4.95)&$\bf1(0.01)$&278(37.76)&209(16.89)&111(5.42)&-(-)&-(-)/182(21.30)&26(0.56)/134(7.47) \\
\hline
$10^{-10}$&$2^{-2}$&8(0.01)&14(0.03)&2(0.00)&$\bf1(0.00)$&7(0.02)&2(0.00)&3(0.00)&7(0.01)&7(0.01)/7(0.01)&$\bf1(0.00)$/3(0.00) \\
&$2^{-3}$&23(0.05)&60(0.16)&2(0.00)&$\bf1(0.00)$&32(0.09)&4(0.01)&3(0.00)&70(0.20)&98(0.27)/46(0.11)&$\bf1(0.00)$/3(0.00) \\
&$2^{-4}$&23(0.06)&156(0.66)&2(0.00)&$\bf1(0.00)$&65(0.24)&4(0.01)&4(0.01)&370(2.31)&437(2.62)/114(0.36)&$\bf1(0.00)$/4(0.01) \\
&$2^{-5}$&16(0.07)&330(4.00)&2(0.00)&$\bf1(0.00)$&122(0.92)&6(0.01)&5(0.01)&-(-)&-(-)/231(2.21)&2(0.01)/5(0.01) \\
&$2^{-6}$&7(0.13)&-(-)&5(0.03)&$\bf1(0.00)$&275(9.61)&15(0.09)&11(0.06)&-(-)&-(-)/313(12.88)&2(0.01)/11(0.06) \\
&$2^{-7}$&4(0.30)&-(-)&20(0.38)&$\bf1(0.01)$&453(84.02)&66(2.26)&38(0.97)&-(-)&-(-)/354(60.34)&7(0.11)/42(1.12) \\
\hline

\end{tabular}
}
\end{center}

\end{table}


\section{Conclusion}\label{Sec4}

We have presented a preconditioner  to  speed up the convergence of the GMRES (FGMRES) iterative method  for solving the system of linear equations arisen from  discretize-then-optimize approach applied to optimal control problems constrained by a Poisson equation. We have shown that the eigenvalues of the preconditioned matrix are one or are located in an interval that the upper bound of the interval is independent of the mesh size. We have also seen that $2\beta$ is lower bound for the non-unit eigenvalues of the preconditioned matrix.  We have compared the numerical results of the proposed preconditioner with those of the several preconditioner presented in the literature. Numerical results show that our preconditioner can be considered as a suitable preconditioner for the investigated problem.

\section*{Acknowledgements}

The authors would like to thank the anonymous referees for their valuable comments and constructive suggestions. The work of Davod Khojasteh Salkuyeh is partially supported by University of Guilan.

\end{document}